\documentclass[12pt]{article}
\usepackage{a4wide}
\usepackage{amsmath,amssymb,amsthm}
\usepackage[T2A]{fontenc}
\usepackage{mathtools}
\usepackage{epsfig}
\usepackage{psfrag}

\author {Ievgen Bondarenko}
\title{\textbf{The word problem in Hanoi Towers groups}}

\newcommand{\HA}{\mathsf{H}}
\newcommand{\HG}{\mathcal{H}}
\newcommand{\A}{\mathsf{A}}
\newcommand{\AG}{\mathcal{G}}

\newtheorem{theorem}{Theorem}

\newtheorem{corollary}{Corollary}[theorem]

\theoremstyle{definition}

\begin{document}
\maketitle

\begin{flushright}
\end{flushright}
%

\begin{abstract}
We prove that elements of the Hanoi Towers groups $\HG_m$ have depth bounded from above
by a poly-logarithmic function $O(\log^{m-2} n)$, where $n$ is the length of an
element. Therefore the word problem in groups $\HG_m$ is solvable in subexponential
time $\exp(O(\log^{m-2} n))$.


\end{abstract}


\section{Introduction}

We consider deterministic finite-state transducers (Mealy automata) with the same input
and output alphabets. Such automata process words over the alphabet letter by letter:
automaton reads the first letter from the current state, produces the output letter, and
changes its state; the output letter and the new state depend only on the current state
and the input letter. In this way, every state $s$ taken as the initial state defines a
transformation of words over the alphabet. If all transformations defined by the states
of an automaton $\A$ are invertible, they generate a group under composition of
functions, which is called the automaton group $\AG_{\A}$ generated by $\A$. Further, we
always assume that the states of automata determine a symmetric generating set of the
group $\AG_{\A}$, so that every element of $\AG_{\A}$ can be given by a word
$s_1s_2\ldots s_n$ over states.

The word problem in every automaton group is solvable. To describe the algorithm, define
the section $s|_v$ of a state $s$ at a word $v$ over the alphabet as the end state of
the automaton after processing the word $v$ from the initial state $s$. This notion
naturally extends to words over the states of an automaton: the section of
$w=s_1s_2\ldots s_n$ at a word $v$ is defined by
\begin{equation}\label{eqn_section_word}
(s_1s_2\ldots s_n)|_v=s'_1s'_2\ldots s'_n, \ \mbox{where } s'_i=s_i|_{(s_{i+1}\ldots
s_n)(v)}
\end{equation}
(we are using left actions). The section $(s_1s_2\ldots s_n)|_v$ has a natural
interpretation: take $n$ copies of the automaton $\A$, choose the initial state $s_i$
in the $i$-th copy, and connect the output of $(i+1)$-th automaton to the input of
$i$-th automaton. Then the final configuration of states after processing a word $v$ is
the section $(s_1s_2\ldots s_n)|_v$. Note that we treat sections as words over states,
together with the natural action on words over the alphabet.

Now the algorithm solving the word problem in automaton groups follows from the fact: a
word $w$ over the states of an automaton $\A$ defines the trivial transformation if and
only if the sections of $w$ at all words over the alphabet act trivially on the
alphabet. If the automaton $\A$ has $k$ states, then $s_1s_2\ldots s_n$ has at most
$k^n$ sections. Therefore this algorithm solves the word problem in at most exponential
time. The precise complexity of the word problem in the class of automaton groups is
unknown. For contracting automaton groups the word problem is solvable in polynomial
time \cite[Proposition~2.13.10]{nekbook}. The word problem in groups generated by
polynomial automata is solvable in subexponential time \cite{bond:poly}. And of course,
if the automaton group is free, nilpotent, etc., then the complexity of the word
problem is smaller than exponential as well.

The complexity of the above algorithm directly depends on the number of sections that
elements have. Define the \textit{section growth function} of an automaton $\A$:
\[
\theta_\A(n)=\max\{ \#Sections(s_1s_2\ldots s_k) : s_i\in \A, k\leq n \}, \ n=1,2,\ldots.
\]
If the function $\theta_\A$ is bounded from above by a polynomial (subexponential)
function then the word problem in the group $\AG_\A$ is solvable in polynomial
(subexponential) time.

The number of sections can be bounded by the depth of state words. Let
$d_\A(s_1s_2\ldots s_n)$ be the least integer $d$ with the property that for every word
$v$ over the alphabet there exists a word $u$ of length $\leq d$ such that
$(s_1s_2\ldots s_n)|_v=(s_1s_2\ldots s_n)|_u$. Define the \textit{depth function} of an
automaton $\A$:
\[
d_\A(n)=\max\{ d_\A(s_1s_2\ldots s_k) : s_i\in \A, k\leq n\}, \ n=1,2,\ldots.
\]
If the alphabet has $a$ letters, then $\theta_\A(n)\leq 1+a^{1}+\ldots+a^{d_\A(n)}$.
Therefore, if the depth function is bounded from above by a poly-logarithmic function,
then the word problem is solvable in subexponential time.

Instead of computing sections as words over states, we can compute a few relations in a
given automaton group, and then reduce sections using these relations. This may highly
reduce the number of sections and, as a consequence, the complexity of the algorithm.
For example, in contracting automaton groups, computing certain relations in a finite
time will guarantee the depth function $d(n)=O(\log n)$ and polynomial word problem.


The goal of this note is to estimate the depth function of the Hanoi automata $\HA_m$
and the complexity of the word problem in the Hanoi Towers groups $\HG_m=\AG_{\HA_m}$,
which model the Tower of Hanoi game on $m$ pegs \cite{gri_sunik:hanoi}. This classical
game is played with $k$ disks of distinct size placed on $m$ pegs, $m\geq 3$.
Initially, all disks are placed on the first peg according to their size so that the
smallest disk is at the top, and the largest disk is at the bottom. A player can move
only one top disk at a time from one peg to another peg, and can never place a bigger
disk over a smaller disk. The goal of the game is to transfer the disks from the first
peg to another peg. For more information about this game, its history, solutions, and
open problems, we refer the reader to \cite{hinz:hanoi,frame-stewart,szegedy} and the
references therein. The automaton model $\HA_m$ presented in \cite{gri_sunik:hanoi}
encodes configurations of disks on pegs by words over the alphabet $\{1,2,\ldots,m\}$,
and the action of each state of the automaton $\HA_m$ corresponds to a single disk move
between two pegs. Therefore each strategy of the player can be encoded by a word over
states of the automaton $\HA_m$. The Hanoi Towers game has subexponential complexity
for $m\geq 4$: it can be solved in $n=exp(O(k^{\frac{1}{m-2}}))$ moves, and this is an
asymptotically optimal solution (see more precise estimates in
\cite{frame-stewart,szegedy}). Note that if we express the height of the tower $k$ in
terms of the length $n$ of an optimal solution, we get $k=O(\log^{m-2} n)$. I do not
see how the estimate on the game's complexity immediately implies the estimate on the
complexity of the word problem in the groups $\HG_m$. Nevertheless, we prove the
following results.

\begin{theorem}\label{thm_main}
The depth function of the Hanoi automaton $\HA_m$ satisfies $d(n)=O(\log^{m-2} n)$.
\end{theorem}

\begin{corollary}
The section growth function of the Hanoi automaton $\HA_m$ satisfies
$\theta(n)=\exp(O(\log^{m-2} n))$.
\end{corollary}

\begin{corollary}
The word problem in the Hanoi Towers group $\HG_m$ is solvable in subexponential time
$\exp(O(\log^{m-2} n))$.
\end{corollary}

\section{Automaton groups and Hanoi automata}

In this section we briefly review necessary information about automata, automaton
groups, and describe the construction of Hanoi automata $\HA_m$. See
\cite{GNS:ADG,gri_sunik:branching,nekbook} for more details.

Let $X$ be a finite alphabet and $X^{*}$ be the free monoid over $X$. The elements of
$X^{*}$ are finite words $x_1x_2\ldots x_n$, $x_i\in X$, $n\in\mathbb{N}\cup\{0\}$, the
identity element is the empty word $\emptyset$, and the operation is concatenation of
words. The length of $v=x_1x_2\ldots x_n$ is $|v|=n$.

An automaton $\A$ over the alphabet $X$ is a finite directed labeled graph, whose
vertices are called the states of the automaton, and for each vertex $s\in \A$ and
every letter $x\in X$ there exists a unique outgoing arrow at $s$ labeled by $x|y$ for
some $y\in X$.

Every state $s\in \A$ defines a transformation of $X^{*}$ as follows. Given a word
$v=x_1x_2\ldots x_n\in X^{*}$, there exists a unique directed path in the automaton
$\A$ starting at the state $s$ and labeled by $x_1|y_1$, $x_2|y_2$,\ldots,$x_n|y_n$ for
some $y_i\in X$. Then the word $y_1y_2\ldots y_n$ is called the \textit{image of
$x_1x_2\ldots x_n$ under $s$}, and the end vertex of this path is called the
\textit{section of $s$ at $v$} denoted $s|_v$. A word $s_1s_2\ldots s_n$ over states
acts on $X^{*}$ by composition: $(s_1s_2\ldots s_n)(v)=(s_1s_2\ldots s_{n-1})(s_n(v))$.
The section of a word $s_1s_2\ldots s_n$ over states at a word $v\in X^{*}$ is defined
by Equation (\ref{eqn_section_word}). Further, by \textit{states in section}
$(s_1s_2\ldots s_n)|_v$ we mean states $s'_i$ given by Equation
(\ref{eqn_section_word}).

If all transformations defined by the states of $\A$ are invertible, the automaton $\A$
is called invertible, and the group $\AG_\A$ generated by these transformations under
composition of functions is called the automaton group generated by $\A$.

\begin{figure}
\begin{center}
\psfrag{a12}{$a_{(12)}$} \psfrag{a13}{$a_{(13)}$} \psfrag{a14}{$a_{(14)}$}
\psfrag{a23}{$a_{(23)}$} \psfrag{a24}{$a_{(24)}$} \psfrag{a34}{$a_{(34)}$}
\psfrag{id}{$e$} \psfrag{1|1}{$1|1$} \psfrag{2|2}{$2|2$} \psfrag{3|3}{$3|3$}
\psfrag{4|4}{$4|4$} \psfrag{1|2}{$1|2$} \psfrag{2|1}{$2|1$} \psfrag{1|3}{$1|3$}
\psfrag{3|1}{$3|1$} \psfrag{1|4}{$1|4$} \psfrag{4|1}{$4|1$} \psfrag{2|3}{$2|3$}
\psfrag{3|2}{$3|2$} \psfrag{2|4}{$2|4$} \psfrag{4|2}{$4|2$} \psfrag{3|4}{$3|4$}
\psfrag{4|3}{$4|3$}\epsfig{file=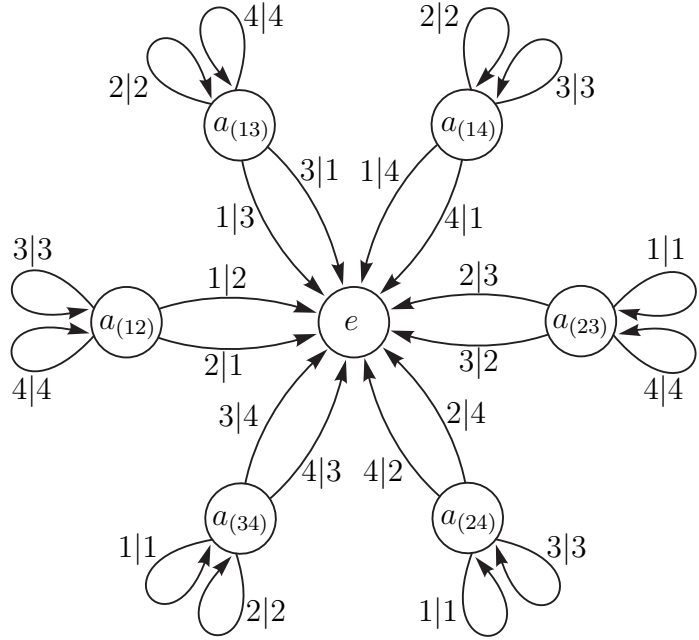,width=260pt} \caption{The Hanoi
automaton $\HA_4$}\label{fig_AutomHanoi}
\end{center}
\end{figure}

The \textit{Hanoi automaton $\HA_m$} is defined over the alphabet
$X_m=\{1,2,\ldots,m\}$. It has the trivial state $e$ and the state $a_{(ij)}$ for every
transposition $(i,j)$ on $X$. All arrows outgoing from $e$ end in $e$ and are labeled
by $x|x$ for each $x\in X$. Each state $a_{(ij)}$ has two outgoing arrows
$a_{(ij)}\rightarrow e$ labeled by $i|j$ and $j|i$, and the other arrows are loops at
$a_{(ij)}$ labeled by $x|x$ for every $x\in X\setminus\{i,j\}$. For example, the
automaton $\HA_4$ is shown in Figure~\ref{fig_AutomHanoi} (the loops at the trivial
state $e$ are not drawn).

The automaton $\HA_m$ is invertible and produces a symmetric generating set of the
Hanoi Towers group $\HG_m=\AG_{\HA_m}$. The action of states $a_{(ij)}$ on the space
$X_m^{*}$ can be given recursively as follows:
\[
a_{(ij)}(iv)=jv, \quad a_{(ij)}(jv)=iv, \quad a_{(ij)}(xv)=xa_{(ij)}(v) \ \mbox{ for
$x\not\in\{i,j\}$}.
\]
In other words, the state $a_{(ij)}$ changes the first occurrence of letter $i$ or $j$
to the other one, and leaves the other letters unchanged. The states of the Hanoi
automata satisfy the following important property: for any $x\in X_m$ and $s\in \HA_m$,
\begin{equation}\label{eqn_property_states}
\mbox{ if $s(x)\neq x$ then $s|_x=e$},\quad \mbox{ if $s(x)=x$ then $s|_x=s$}.
\end{equation}

\section{Proof of Theorem~\ref{thm_main}}

\begin{proof}
Let $d_m$ be the depth function of the Hanoi automaton $\HA_m$. We will prove
\begin{equation}\label{eqn_main_estimate}
d_{m}(n)\leq m^2(m-1)^2\ldots 3^2 (\log n)^{m-2}+d_{m-1}(n), \qquad d_3(n)\leq \log n+1.
\end{equation}
(For simplicity, to avoid extra brackets, by $\log n$ we mean the least integer greater
than the binary logarithm of $n$).

\textbf{Case $m=3$.} Let us estimate the depth of certain words over states. If
$w=ee\ldots e$ then $d_3(w)=0$. We say that a word $w=s_1s_2\ldots s_n$ is a one-letter
word, if there exists $s\in \HA_3$ such that $s_i\in\{e,s\}$ for all $i$. Note that
$d_3(w)\leq 1$ for every one-letter word $w$. If $w_1, w_2$ are one-letter words (maybe
for different letters $s$), then the properties (\ref{eqn_property_states}) imply that
the section $(w_1w_2)|_x$ is a one-letter word for every $x\in X_3$, and
$d_3(w_1w_2)\leq 2$. It follows that if a word $w$ is a concatenation of $n$ one-letter
words, then $w|_x$ is a concatenation of at most $(n+1)/2$ one-letter words. The
estimate $d_3(n)\leq \log n+1$ follows.

\textbf{Case $m>3$.} Fix $x\in X_m$. Let $\HA_m^{(x)}$ be the set of all states $s\in
\HA_m$ that fix $x$. Note that $\HA_m^{(x)}$ is a subautomaton of $\HA_m$. In
particular, for any $s_i\in \HA_m^{(x)}$ and $v\in X_m^{*}$ every state in the section
$(s_1s_2\ldots s_n)|_v$ fixes $x$. Application of inductive hypothesis is based on the
following observation: If we restrict the automaton $\HA_m^{(x)}$ to the alphabet
$X_m\setminus\{x\}$, we get the automaton $\HA_{m-1}$. This gives an estimate on the
number of sections at words over $X_m\setminus\{x\}$. But since $s(x)=x$ and $s|_x=s$
for each $s\in \HA_m^{(x)}$, we have
\[
(s_1s_2\ldots s_n)|_v=(s_1s_2\ldots s_n)|_u
\]
for every $v\in X_m^{*}$ and $s_i\in \HA_m^{(x)}$, where $u$ is the word over
$X_m\setminus\{x\}$ that is obtained from $v$ by removing every occurrence of letter
$x$. Therefore it is sufficient to count only sections at words over
$X_m\setminus\{x\}$. We have proved the estimate
\begin{equation}
d_{m}(s_1s_2\ldots s_n)\leq d_{m-1}(n)
\end{equation}
for all $s_i\in \HA_m^{(x)}$, $n\in\mathbb{N}$, $m>3$.

Now estimate (\ref{eqn_main_estimate}) follows once we prove the following statement.

We say that a \textit{section $(s_1s_2\ldots s_n)|_v$ satisfies the property (*)} if
there exists $x\in X_m$ such that every state in $(s_1s_2\ldots s_n)|_v$ fixes $x$,
i.e., belongs to $\HA_m^{(x)}$.

\textbf{Claim.} For every $n\in\mathbb{N}$ and any states $s_1,\ldots,s_n\in \HA_m$, the
section $(s_1s_2\ldots s_n)|_v$ satisfies the property (*) for all words $v\in X_m^{*}$
of length $|v|\geq C_m(\log n)^{m-2}$ with constant $C_m=3^24^2\ldots m^2$. In
particular, $(s_1s_2\ldots s_n)|_v$ is a word over the states of $\HA_{m}^{(x)}$ for
some $x\in X$, and $d_{m}((s_1s_2\ldots s_n)|_v)\leq d_{m-1}(n)$.

\begin{proof}

We prove the claim by induction on $m$. The basis of induction $m=3$ is shown above.
Suppose the claim holds for less than $m$ and consider the case $m$. We are going to use
the inductive hypothesis as follows: for every $x\in X_m$, any $s_i\in \HA_m^{(x)}$, and
all words $v\in X^{*}_m$ that contain at least $C_{m-1}(\log n)^{m-3}$ letters different
from $x$, there exists a letter $y\in X_m$, $y\neq x$, such that every state in the
section $(s_1s_2\ldots s_n)|_v$ fixes $y$.

Take two different letters $x,y\in X_m$ and elements $g=s_1s_2\ldots s_k$ for $s_i\in
\HA_m^{(x)}$ and $h=t_1t_2\ldots t_l$ for $t_i\in \HA_m^{(y)}$. Consider the section
$(gh)|_v$ for a word $v\in X^{*}$ of length $\geq |X_m|C_{m-1}(\log n)^{m-3}$, $n=k+l$.
Then $v$ contains at least $C_{m-1}(\log n)^{m-3}$ letters $z$ for certain $z\in X_m$.
If $z\neq y$ then we can apply inductive hypothesis to $h|_v$: every state in $h|_v$
fixes some letter besides $y$. If $z=y$ then $h(v)$ contains at least $C_{m-1}(\log
n)^{m-3}$ letters $z=y\neq x$, and every state in $g|_{h(v)}$ fixes some letter besides
$x$. If we get a common letter for $g|_{h(v)}$ and $h|_v$, then $(gh)|_v$ satisfies the
property (*). Otherwise, we get at least three letters, each one fixed either by all
states in $h|_v$, or in $g|_{h(v)}$. We can proceed further and consider $(gh)|_v|_u$
for words $u$ of length $\geq |X_m|C_{m-1}(\log n)^{m-3}$. Then either $(gh)|_v|_u$
satisfies the property (*), or there are at least four letters, each one fixed either
by all states in $h|_v$, or in $g|_{h(v)}$. It follows that for all words $v$ of length
$\geq |X_m|^2C_{m-1}(\log n)^{m-3}$ the section $(gh)|_v$ satisfies the property (*).

Now consider any word $w=s_1s_2\ldots s_n$, $s_i\in \HA_m$. We partition $w$ on subwords
\begin{equation}\label{eqn_blocks}
w=w_1w_2\ldots w_k, \quad k\leq n
\end{equation}
such that every subword $w_i=s_{j_i}s_{j_i+1}\ldots s_{j_{i+1}-1}$ satisfies the
property (*). Consider the words $w_1w_2$, $w_3w_4$, \ldots, and their sections at words
$v$ of length $\geq |X_m|^2C_{m-1}(\log n)^{m-3}$. Then, using the fact proved in the
previous paragraph, the section $w|_v$ can be represented in the form~(\ref{eqn_blocks})
with $\leq (k+1)/2$ subwords. Applying this procedure $\log k\leq\log n$ times, we get
that $g|_v$ satisfies the property (*) for all words $v$ of length $\geq
|X_m|^2C_{m-1}(\log n)^{m-2}$. The claim is proved.
\end{proof}
\end{proof}

Using program package \cite{AutomGrp} we have calculated the values of the depth
function $d_4(n)$ and the section growth function $\theta_4(n)$ of the Hanoi automaton
$\HA_4$ for small values of $n$:
\begin{center}
\begin{tabular}{|c|c|c|c|c|c|c|c|c|c|c|c|c|}
  \hline
  $n$ & 1 & 2 & 3 & 4 & 5 & 6 & 7 & 8 & 9 & 10 & 11 & 12 \\
  \hline
  $d_4(n)$ & 1 & 2 & 2 & 3 & 4 & 4 & 5 & 5 & 6 & 6 & 6 & 7 \\
  \hline
  $\theta_4(n)$ & 2 & 4 & 8 & 13 & 17 & 24 & 31 & 39 & 48  & 60 & 70 & 81\\
  \hline
\end{tabular}
\end{center}
This suggests that $\theta_4(n)$ may have polynomial (quadratic?) growth, which implies
polynomial word problem.


\end{document}